
\magnification=\magstep1
\overfullrule=0pt
\def\n{\noindent}
\def\pf {\n{\bf Proof.} }
\def\prm{\prime}
\def\ss{\smallskip}
\def\ms{\medskip}
\def\msk{\medskip}
\def\bs{\bigskip}
\def\bsk{\bigskip}

\def\scon{\hbox{${\cal C}_n$}}

\def\R{\mathop{{\rm I}\kern-.2em{\rm R}}\nolimits }
\def\E{\mathop{{\rm I}\kern-.2em{\rm E}}\nolimits }
\def\Rn{\hbox{${\mathop{{\rm I}\kern-.2em{\rm R}}\nolimits}^{\hbox{\rm n}}$} }

\def\P{\mathop{{\rm I}\kern-.2em{\rm P}}\nolimits}
\def\d{\, d}
\def\pn{\par\noindent}
\def\prod{\mathop{\Pi}}

\def\bbbox{\quad{\vrule height4pt depth0ptwidth4pt}}
\def\dsp{\displaystyle}
\def\tF{\tilde F}
\def\Bn{{\rm B_2^n}}
\def\Sn{{{\rm S}^{n-1}}}
\def\On{{{\rm O}(n)}}
\def\N{{\rm I\kern-1.6pt{\rm N}}}
\centerline{\bf
{{ On\ the\ Gaussian\ measure\ of\ the\ intersection}}}
\centerline{\bf
{{ \ of\ symmetric,\ convex\ sets}}}
\bs
\centerline {by}
\bs
\centerline{G. Schechtman\footnote{$^1$}{Supported  
in part by US-Israel Binational Science Foundation}${^,} 
$\footnote{$^{2}$}{Participant, NSF workshop, Texas A\&M University}, 
 Th. Schlumprecht \footnote{$^{3}$}{Supported in part by NSF Grants 
 DMS-9203753    
 and DMS-9501243}${^,}$\footnote{$^4$}{Supported in part by Texas Advanced 
 Research Program Grant 160766} and J. Zinn$^{1,4,}$\footnote{$^5$}{Supported in part by NSF Grant DMS-9208053}}
\footnote{}{{\sl AMS} 1991 {\it Subject Classification.} Primary 60E15, 28C20}
\ms
\centerline{\it Weizmann Institute of Science, Texas A\&M University, Texas A\&M University}
\bigskip
\par
{\narrower\narrower \noindent \sevenrm\ \ \ \ \  The Gaussian Correlation Conjecture states that for any two symmetric,
convex sets in n-dimensional space and for any centered, Gaussian
measure on that space, the measure of the intersection is greater than
or equal to the product of the measures. In this paper we obtain
several results which substantiate this conjecture. For example, in the standard Gaussian case, we show there is a positive constant, c, 
such that the conjecture is true if the two sets are in the 
Euclidean ball of radius $c\sqrt{n}$. Further we show that 
if for every n the conjecture is true when the sets are in
the Euclidean ball of radius $\sqrt{n}$, then it is true in general.
Our most concrete result is that the conjecture is true if the two
sets are (arbitrary) centered ellipsoids.\par}

\bs
\n
{\bf Introduction.} The standard Gaussian measure on $\Rn$ is given by its
density:
$$\mu_n(A)={1\over (2\pi)^{n/2}}\int_A e^{-|x|^2/2} \d x.$$
\pn A general mean zero Gaussian measure on \Rn\ is  a linear image of
the standard Gaussian measure.
\par
Let $\scon$ denote the collection of convex closed subsets of \Rn which are
 symmetric about the origin.

\proclaim Conjecture C. For any $n\ge 1$, if $\mu$ is a mean zero, Gaussian
 measure on $\Rn$, then for all $A, B\in \scon$,
$$\mu(A\cap B)\ge\mu(A)\mu(B).$$\par

Recall that a function $f:\Rn\to \R^+$ is called {\it quasi concave} if for
  any $r\in\R$ the set $\{x\in\Rn: f(x)\ge r\}$ is convex. For such an $f$ let 
$A=\{(x,t):f(x)\ge t\}$ and $A_t=\{x: f(x)\ge t\}$. Then, $A_t$ is convex and symmetric if $f$ is symmetric and further, 
$$f(x)=\int_0^\infty I_{A_t}(x)\, dt.$$
By Fubini's theorem  Conjecture (C) has the following
functional version.

\proclaim Conjecture C$^\prm$. Let $f,g$ be  non-negative, quasi-concave, and symmetric. Then
 $$\E_{\mu_n}(f\cdot g)\ge \E_{\mu_n}(f)\cdot\E_{\mu_n}(g),$$
where  $\E_{\mu_n}(f)$ denotes the expectation of $f$ with respect to $\mu_n$.

\n  It is, of course, enough to show conjecture (C) for symmetric and convex 
    polytopes. Since convex, symmetric polytopes are images of the 
    unit cube $[-1,1]^m$ in some possibly higher dimensional space, $\R^m$, 
    under a linear map an easy integral transformation shows that (C) is
    equivalent to the following conjecture (C$^{\prm\prm}$) which is stated in a more probabilistic language.
\proclaim Conjecture C$^{\prm\prm}$. 
 If $\{X_i\}_{i=1}^n$ are jointly Gaussian, mean zero random variables,
 and \hfil\break
$1\le k\le n$ then,
$$P(\max_{i\le n}|X_i|\le 1)\ge
 P(\max_{i\le k}|X_i|\le 1)P(\max_{k<i\le n}|X_i|\le 1).$$\par

According to Das Gupta, Eaton, Olkin, Perlman, Savage and Sobel [DEOPSS],
the history of this
problem prior to 1970 starts with a paper of Dunnett and Sobel [DS] in
(1955) and after contributions by Dunn [Du] in (1958), it culminated in
papers of Khatri [Kh] and \u Sid\' ak [Si1], both in (1967), in which they
independently obtained (C$^{\prm\prm}$) in the case $k=1$:

\par\n{\bf Theorem}  (Khatri, \u Sid\' ak).{\sl  Let $\{X_i\}_{i=1}^n$ be jointly 
 Gaussian, mean zero random variables.  Then
$$P(\max_{i\le n}|X_i|\le 1)\ge 
  P(|X_1|\le 1)P(\max_{1<i\le n}|X_i|\le 1).$$\par}
\msk
\n If a symmetric slab is defined to be a set of the form $\{x\in \Rn:
 |<x,u>|\le 1\}$ for some $u\in \Rn$, Theorem 1 is equivalent to
\msk
\proclaim Theorem. If $\mu$ is a mean zero Gaussian measure on \Rn
 , $A \in\scon$ , and $S$ is a symmetric slab, then
$$\mu(A\cap S)\ge\mu(A)\mu(S).$$\par

\n As a corollary of the theorems above, they obtained a result which
solved the problem
 studied by Dunnett and Sobel [DS] and  Dunn [Du].
\msk
\n{\bf Corollary}  (Khatri, \u Sid\' ak).{\sl 
  $$P(\max_{i\le n}|X_i|\le 1)\ge\prod_{i=1}^n
 P(|X_i|\le 1).$$\par}
\msk
Another important milestone for this problem was achieved by 
 the work of L.~D.~Pitt in
 1977, where the two dimensional case was settled.
\msk
\par\n{\bf Theorem} ([Pi]).{\sl For any $A,B\in{\cal C}_2$   
   $\mu_2(A\cap B)\ge\mu_2(A)\mu_2(B)$ .}
\msk

In [DEOPSS] and  Gluskin [Gl] measures other than Gaussian measures are
considered. The problem can and has been attacked using 
measure theoretic, geometric and analytic techniques. In this
note we present several   partial results using some of these
techniques. In particular we prove, in section 1, that the
conjecture is true for ``small enough" sets. We also show, in Proposition 5, that the result
holds ``on the average". 
It follows from the remark following Proposition 1 that, if, in the
statement of conjecture C, one puts the factor $2^{n/2}$ on the left hand
side,  then the resulting statement is true. In Proposition 7 we prove that
if one could replace the factor
$2^{n/2}$ with
$2^{o(n)}$, then the conjecture would follow.
In Theorem 8 (section 2) we prove the 
conjecture for sets more general than sets having a common ``orthogonal
unconditional" basis. 
Finally, in Theorem 10 we show that the conjecture holds for arbitrary  
centered ellipsoids in $\Rn$.  

We will need the following notations and concepts. In $\Rn$ the usual unit basis will be denoted by 
 $e_1,e_2,...e_n$ , $|\cdot|$ is
  the Euclidean norm, and $<\cdot,\cdot>$ the scalar product generated by 
  $|\cdot|$. $\Bn=\{x\in\Rn:|x|\le 1\}$ will be the Euclidean unit ball
   and $\Sn=\{x\in\Rn:|x|\le 1\}$ its sphere. The orthogonal group on $\Rn$,
   i.e. the set of real unitary $n\times n$ matrices, will be denoted by
   $\On$.  Lebesgue measure
   on $\Rn$ will be denoted by $m_n$.
   
   We will make heavy use of the following concept from convex geometry.
Recall that a non-negative function
$f:\Rn\to\R^+$ is called
 log-concave if  for  $x,y\in\Rn$ and  $0\le t\le 1$,
$$ f(tx+(1-t)y)\ge  f(x)^tf(y)^{1-t},$$ i.e $\log f$ is concave on its 
  support.

 Note that the indicator functions of convex sets are log-concave and
  that log-concave functions are quasi-concave. We also will 
 need  the following deep result of Pr\'ekopa and Leindler.

\par\n{\bf Theorem}([Le] and [Pr], see also  [BL]). 
{\sl If $f$ is log-concave on
\Rn
 and
$1\le k<n$, then the function $g:{\R}^{k}\to \R$ , with
$$g(x_1,\ldots, x_k)=\int_{{\R}^{n-k}}
f(x_1,\ldots,x_k,z_1,\ldots,z_{n-k}) \d z$$ is also log concave.}
\pn

 Since $h\circ A$ is log concave whenever $h$ is log concave and $A$ is linear,
  and since the product of two log concave functions is also log concave
  the following Corollary follows immediately. 
\medskip
\proclaim Corollary. If $f$ and $g$ are log concave, so is 
 $y\mapsto \int f(x+y)g(x)\d x$.\par
\ss In order to get a glimpse of the mysterious power of the 
  Pr\'ekopa-Leindler result  we will use it in order to  
  give a  very short proof of the result of Khatri and \u Sid\' ak.

  We first observe that the conjecture (C) and  thus (C$^\prm$) are trivially true
   in the case $n=1$.
  Assume that $S=\{x\in\Rn: |x_1|\le s\}$ and that $A\in\scon$. For $x_1\in\R$,
  
  $f(x_1):=\int_{R_{n-1}} I_A(x_1,y) \d\mu_{n-1}(y)$.
   Since the density of $\mu_{n-1}$ and $I_A$ are log concave we
    deduce from [Le] and [Pr] that $f$ is a log concave function  
    on $\R$ and thus
    $$\mu(A\cap S)=\int_{\R} I_{[-s,s]} f(x_1)\d\mu_1( x_1)\ge
     \mu_1([-s,s])\cdot \E_{\mu_1}(f)=\mu(S)\cdot\mu(A),$$
     where the inequality follows from the one dimensional case.

\bsk

\n{\bf Section 1. Restriction on size.} 
\msk
Using the rotation on 
$\Rn\times\Rn$ given by $(x,y)\mapsto({x+y\over\sqrt{2}},{x-y\over\sqrt{2}})$ leads to
the following observation.

\proclaim  Proposition 1. If $A, B\in\scon$, we have
$$\mu_n(A)\cdot\mu_n(B)\le\mu_n(\sqrt{2}(A\cap B))\mu_n({(A+B)\over
\sqrt{2}}).$$

\par\n
{\bf Proof.} Using the rotational invariance of the measure $\mu_n\otimes\mu_n$ we get

$$\eqalign{\mu_{2n}(A\times B)&=\int I_A(x)I_B(y)\d\mu_n(x)\d 
\mu_n(y)\cr &=\int I_A({u+v\over \sqrt{2}})I_B({v-u\over \sqrt{2}})\,
\mu_n(du) \mu_n(dv)
 \cr &=\int \mu_n((\sqrt{2}A-u)\cap(\sqrt{2}B+u))\, \mu_n(du).}$$
 \par\n
 Note that for $u\in\Rn$ it follows that
 $(\sqrt{2}A-u)\cap(\sqrt{2}B+u)$ is not empty if and only if
 there exists  a $z\in\Rn$ for which ${z+u\over \sqrt{2}}\in A$ and
 ${z-u\over \sqrt{2}}\in B$. Since that can only happen if
  $u$ lies in $(A-B)/\sqrt{2}=(A+B)/\sqrt{2}$ we deduce
  that the integrand can only be non zero on ${(A+B)/\sqrt{2}}$.
\par\n
Furthermore, the mapping 
$u\mapsto\int \mu_n((\sqrt{2}A-u)\cap(\sqrt{2}B+u))\,   
 \mu_n(du)$ is log concave by the Pr\'ekopa-Leindler theorem. Since  
  it is also symmetric, it  
 is maximized at zero. 
 Hence the integral is bounded by
 $\mu_n(\sqrt{2}(A\cap B))\cdot\mu_n({(A+B)\over \sqrt{2}}).\bbbox$
\msk
\par\n{\bf Remark}.
 Note that for any measurable $K\subset\Rn$ and $c>1$ it follows that
  $\mu_n(cK)=(2\pi)^{-n/2}\int I_K(x/c)\cdot e^{-|x|^2/2}\d x=
  c^n(2\pi)^{-n/2}\int I_K(u)\cdot e^{-c^2|u|^2/2}\d u\le c^n\mu_n(K)$. 
  Thus Proposition 1 implies  
$\mu_n(A)\mu_n(B)\le 2^{n/2}\mu_n(A\cap B)$ if $A,B\in\scon$.

  Using 
   $m_n(\cdot)\ge (2\pi)^{-n/2}\mu_n(\cdot)$, we deduce the following corollaries.
\msk
\proclaim Corollary 2. For
 $A,B\in\scon$  we have
$$\mu_n(A\cap B)\ge {(2\pi)^{n/2}\over m_n(A+B)}\mu_n(A)\mu_n(B).$$\par
\ms

\proclaim Corollary 3. Suppose $\rho_n$ is chosen so that 
$m(2\rho_n \Bn)=(2\pi)^{n/2}$. (Note that
$\rho_n={\dsp {1\over \sqrt{2}}}(\Gamma(1+{\dsp {n\over
2}})^{1/n}\sim{\dsp {1\over 2} \sqrt{n\over e}}$.)
 \par\n
 {\sl Then, $\mu(A\cap B)\ge\mu(A)\mu(B)$, for all $A,B\in\scon$ with
  $A,B\subset\rho_n\Bn$.}
\ms
In Corollary 7 below we will show that, if we could replace the factor $\rho_n$ by
$\sqrt{n}$, then the conjecture would follow. We first make the following 
observation which indicates that it would be enough to show 
 (C) approximately.

\proclaim  Proposition 4. Assume that there is a sequence of positive numbers $(c_n)$ with \hfil\break
$\lim_{n\to\infty} c_n^{1/n}=1$, so that
$\mu_n(A\cap B)\ge c_n\mu_n(A)\mu_n(B)$,
  for all $n\in\N$ and $A, B\in\scon$. 
  Then, for all $n\in\N$ and $A, 
  B\in\scon$,
$$\mu_n(A\cap B)\ge \mu_n(A)\mu_n(B).$$\par

\n {\bf Proof.} For each $N$ consider $A^N=A\times\cdots\times A$, and
$B^N$. The
 assumption gives:
$$
\mu_n^N(A\cap B)=\mu_{Nn}(A^N\cap B^N)\ge c_{Nn}\,\mu_{n}(A)\mu_n(B).
$$
Taking $N^{\hbox{th}}$ roots, letting $N\to\infty$ and using the
hypothesis, the
 result follows.\bbbox

We    now  show that the conjecture holds on the average. This is true for 
 more general measures and more general sets.

\proclaim Proposition 5. Let $m$ be the Haar measure on the orthogonal
group $\On$, and let $\nu$ be a rotational invariant probability on $\Rn$
 assume that $A,B\subset\Rn$ are two star shaped sets with 0 being a center,
  i.e. for any $\theta\in\Sn$ the set $\{r\ge 0: r\theta\in A\}$
   is an interval, which we will denote by $A_\theta$.
 \par\n 
  {\sl  Then it follows that 
  $$   \int_\On \nu(A\cap U(B)) \d m(U)\ge \nu(A)\nu(B).$$}
  
  \ss\n
  {\bf Proof.} Since $\nu$ is rotational invariant it is the
  image of some product probability $\nu_1\otimes\sigma_n$ ($\nu_1$ being a 
  probability 
   on $[0,\infty)$) under the map: 
   $\Sn\times[0,\infty)\ni(\theta,r)\mapsto \theta r$. We will also use the 
    fact that for any $\theta_0$ the  measure $\sigma_n$  is the image
     of $m$ under the map $\On\ni U\mapsto U(\theta_0)$. Finally
     we observe that for two star shaped sets $A$ and $B$, with 0 being 
     their center, and for any two $\theta$, and $\theta'$ we deduce that
      $\nu_1(A_\theta\cap B_{\theta'})=\min(\nu_1(A_\theta),\nu( B_{\theta'}))
       \ge \nu_1(A_\theta)\cdot\nu( B_{\theta'})$.
  
  These observations allow us to make the following estimates.
  $$\eqalign { \int_\On \nu(A\cap U(B)) \d m(U)&=
       \int_\Sn\int_\Sn\int_0^\infty 
        I_{A_\theta}(r)I_{B_{\theta'}}\d \nu_1(r)\d\sigma_n(\theta)
	         \d \sigma_n(\theta')   \cr
  &=\int_\Sn\int_\Sn \nu_1(A_\theta\cap B_{\theta'})	\d\sigma_n(\theta)
	         \d \sigma_n(\theta')   \cr
   &\ge\int_\Sn\int_\Sn \nu_1(A_\theta)\nu_1(\cap B_{\theta'})	
                \d \sigma_n(\theta)         \d \sigma_n(\theta')   \cr
   &=\int_\Sn \nu_1(A_\theta)\d \sigma_n(\theta)   
      \int_\Sn \nu_1(B_{\theta'})\d \sigma_n(\theta')=\nu(A)\nu(B),    \cr}$$
      which proves the claim. 
    \bbbox  
  \ss\n 
\proclaim Corollary 6.  For any $r>0$ and any $A\in \scon$,
$$
\mu_n(A\cap r\Bn)\ge \mu_n(A)\mu_n(r\Bn).
$$
\bigskip

\n Here is one example of how to use the above results.

\proclaim Corollary 7. If for all $n$, $\mu_n(A\cap B)\ge\mu_n(A)\mu_n(B)$  for
 all $A, B\in \scon$ for which $A, B\subset\sqrt{n}\Bn$, then the inequality holds for all $n$ and 
 $A,B\in\scon$.\par
\n{\bf  Proof.}  For $A, B\in\scon$, we have 
$$\eqalign{\mu_n(A\cap B)
    \ge&\mu_n(A\cap B\cap \sqrt{n}\Bn)
            \ge\mu_n(A\cap\sqrt{n}\Bn)\mu_n(B\cap\sqrt{n}\Bn)\cr
   \ge&\mu_n(A)\mu_n(B)\mu_n^2(\sqrt{n}\Bn),}$$ 
 by Corollary 6.
From the Central Limit Theorem we deduce,
$$\mu_n(\sqrt{n}\Bn)=\mu_n(\sum_{i=1}^nx_i^2\le n)
  =\mu_n({\dsp {\sum_{i=1}^n(x_i^2-1)\over \sqrt{n}}}\le 0)\to 1/2,$$
so the above Proposition applies with $c_n=\mu_n(\sqrt{n}\Bn)$.\bbbox
\ss\n
{\bf Remark.} In view of Corollary 7 it would seem that one could
improve upon the last result (using a ball of a somewhat smaller radius
than
$\sqrt{n}$). However the improvement is negligible.
\ss

\n{\bf Section 2. Geometrical restrictions.}
\ss
By induction on the dimension  it is
easy to see that the conjecture is true if the convex symmetric sets are
 1-unconditional with respect to the same orthogonal basis
$\{e_i\}_{i=1}^n$
 (i.e., $(x_1,\ldots,x_n)\in A\iff (\pm x_1,\ldots,\pm x_n)\in A)$.   
 Here we relax somewhat the geometrical restrictions.

\proclaim Proposition 8. Let $\nu$ be a product probability on $\Rn$. If $A, B\in\scon$ satisfy: 
{\item{(i)} $x\in A\cap B\Longrightarrow x_ie_i\in
A\cap B, \forall i\le n$
\item{(ii)}for every pair of
orthants,   $Q$ and $Q'$, 
$[\nu(A\cap Q)-\nu(A\cap Q^{\prime})][\nu(B\cap
Q)-\nu(B\cap Q^{\prime})]\ge 0$. (in particular, if $\nu(B\cap Q)$ are
all equal).\hfil\break}
 Then, $\nu(A\cap B)\ge \nu(A)\nu(B)$.

\n To prove this we need the following result. It 
can be found in [KR] and is related to a result in [AD].\par\n
\ms
\n{\bf Theorem.} ([KR] ).
{\sl Let $\nu$ be a  product measure on $\Rn$ and let $f_i, 1\le i\le 4$, be
non-negative functions on
 $\Rn$ satisfying: 
  $$f_1(x)\cdot f_2(y)\le f_3(x\vee y)\cdot f_4(x\wedge y).$$
 Then 
  $$\int f_1 d\nu\cdot\int f_2 d\nu\le \int f_3 d\nu\cdot\int f_4 d\nu.$$}
\ms
\n {\bf Proof of Proposition 8.}   We shall first  prove  that the Karlin-Rinott  theorem implies that, for each orthant
$Q$,
$$
\nu(A \cap Q)\nu( B\cap Q)\le
          \nu(Q)\nu(A\cap B\cap Q). \leqno{(1)}
$$
  Let
$Q$ represent an orthant, say the first
 orthant, and let $f_1=I_{A\cap Q}, f_2=I_{B\cap Q}, f_3=I_Q $ and
$f_4=I_{A\cap
 B\cap Q}$. To use the  Karlin-Rinott theorem we need to show
$$
x\in A\cap Q, y\in B\cap Q\Longrightarrow x\vee y\in Q \hbox{ and } x\wedge
 y\in A\cap B.
$$ Without loss of generality we may assume that $x$ and $y$ are in the
interiors of
$A\cap Q$ and $B\cap Q$, respectively.
 We need to show that $x\wedge y\in A\cap B$. Assuming this were not true
 we let $w$ be the point in $ A\cap
B\cap Q$ which is the closest to $x\wedge y$. By the Pythagorean theorem,
$w_i\le(x\wedge y)_i$ for every $1\le i\le n$. By (i), the rectangular box
$R=\{z\in Q ;\ z_i\le w_i,\  \forall i\}$ is contained in $A\cap B$. Let
$U$ be an open set such that $x\in U\subseteq A\cap Q$ and similarly
$V$  an open set such that $y\in V\subseteq B\cap Q$. $w$ is an interior
point of the convex hull of $U$ and $R$ which is a subset of $A$.
Similarly, $w$ is an interior point of the convex hull of $V$ and $R$
which is a subset of $B$. Hence $w$ is an interior point of $A\cap B\cap
Q$. Therefore, if $x\wedge y$ is not already in
$A\cap B\cap Q$, we reach a contradiction.

The Karlin-Rinott theorem now yields (1). Now apply (ii) in order to deduce that
  $$2^{-n}\sum_{Q,Q'} \nu(A\cap Q)\nu(B \cap Q')\le 
    \sum_{Q} \nu(A\cap Q)\nu(B\cap Q),$$
which implies together with (1) the claim.
 \bbbox
\bigskip
 We now want to show the correlation conjecture for 
  two ellipsoids (in arbitrary position).
\proclaim Theorem 9. If $A$ and $B$ are centered ellipsoids in $\Rn$, then 
$\mu_n(A\cap B)\ge \mu_n(A)\mu_n(B)$.

From Theorem 8 it follows that $\mu_n(E\cap F)\ge \mu_n(E)\mu_n(F)$ if  
 $E$ and $F$ are ellipsoids with the same axis. Using the rotational 
 invariance of $\mu_n$ we would be able to deduce Theorem 9 if
  we could show that for two ellipsoids $E$ and $F$ in the standard position,
   i.e. $E=\{x\in\Rn:\sum_{i=1}^n {x_i^2\over r_i^2}\le 1\}$, and 
   $ F=\{x\in\Rn:\sum_{i=1}^n {x_i^2\over \rho_i^2}\le 1\}$,   the minimum
    of $\mu_n(U(E)\cap F)$  over all $U\in\On$ is 
    attained when $U$ is some row permutation of the  identity. Actually
    this is true for all rotational invariant measures on $\Rn$.

\ms
\n{\bf Theorem 10.} {\sl Let $\nu$ be a rotation invariant measure 
on $\Rn$, and
 let 
 $$E=\{x\in\Rn:\sum_{i=1}^n {x_i^2\over r_i^2}\le 1\},{\rm\ and\ }
   F=\{x\in\Rn:\sum_{i=1}^n {x_i^2\over \rho_i^2}\le 1\}$$
 be two ellipsoids in standard position. Then  the value of
 $\min\{ \nu(U(F)\cap E): U \in\On \}$
 is achieved for some row permutation $P$ of the identity, in particular
  this means that $P(F)$ and $E$ are ellipsoids with the same axis. }   
\ss
\pf Using a standard perturbation argument we can and will make the following assumptions.
 
Instead of considering the minimum of the mapping
$\On\ni U\mapsto\int I_E(U(x))I_F(x)\d\nu(x)$
we let $f:[0,\infty)\to[0,\infty)$ be a continuously differentiable function
with $f'(r)<0$ whenever $r>0$, define for $x\in\Rn$
 $\tF(x)=f(|x|^2_F)$ where $|x|^2_F=\sum_{i=1}^n {x_i^2/ \rho_i^2}$ 
 and we assume that $U_0\in\On$  for which
 
 $$\int I_E(U_0(x))\tF(x)\d\nu(x) 
        =\min_{U\in\On}\int I_E(U(x))\tF(x)\d\nu(x).$$

We also assume that the radii $r_1,r_2,\ldots,r_n$ of $E$ and the radii
 $\rho_1,\rho_2,\ldots,\rho_n$ of $F$ are distinct. Finally, we will assume
 that $\nu$ has a strictly positive density $g(|x|)$ with respect to
  $m_n$. 
 
 In order to deduce the claim we will show that the matrix
 $$U_0^T\circ
 \pmatrix{   r_1^{-2} &   {}  &{}   \cr
              {}      &\ddots  &{}   \cr
	      {}      &   {}  &r_n^{-2}\cr} \circ U_0$$
 is diagonal. Since the values  $r_i^{-2}$ are distinct for $i=1,2\ldots n$ 
 this would imply that $U_0$ must be a row permutation of
  some diagonal matrix $J$ which has only the values $1$  or $-1$ in its
  diagonal. Since $J(G)=G$ for any ellipsoid, we can assume that $J$ is
  the identity. 	      

We start with a variational argument. For $i\not= j$ in $\{1,2,\ldots,n\}$ and 
 $\alpha\in\R$, let $V_{(i,j)}^{(\alpha)}$ be the matrix which acts
 on $\Rn$  in the following way. For $x=(x_1,\ldots,x_n)\in\Rn$ we set
 $$V_{(i,j)}^{(\alpha)}(x):=
      (x_1,\ldots,x_{i-1},x_i\cos\alpha-x_{j}\sin\alpha,x_{i+1},\ldots,
               x_{j-1},x_i\sin\alpha+x_j\cos\alpha,x_{j+1},\ldots,x_n),$$  
i.e. $V_{(i,j)}^{(\alpha)}$ acts on the two dimensional subspace of $\Rn$ spanned by $e_i$ and $e_j$ as a rotation by  $\alpha$,
     and on the orthogonal complement of that subspace, it is the identity. 	       	
 
  Using the minimality of $U_0$ we deduce that 
$$\eqalign{ 0=&{\partial\over\partial \alpha}\left[\int I_E(U_0(x))\tF( V_{(i,j)}^{(\alpha)}(x)) g(|x|) \d x \right]_{\alpha=0}\cr
               =&\int I_E(U_0(x)) f'(|x|_F^2)
                 {\partial\over\partial \alpha}
 \left[{(x_i\cos\alpha-x_j\sin\alpha)^2\over\rho_i^2}+{(x_j\cos\alpha+x_i\sin\alpha)^2\over\rho_j^2}\right]_{\alpha=0}g(|x|)\d x                     \cr
               =&2(\rho_j^{-2}-\rho_i^{-2})\int x_i x_jI_E(U_0(x)) f'(|x|_F^2)
                          g(|x|)\d x.                     \cr}$$   
  We fix $i\le n$, and  for $x=(x_1,\ldots,x_n)\in\Rn$ we let
   $x^{(i)}=(x_1,\ldots,x_{i-1},x_{i+1},\ldots,x_n)\in\R^{n-1}$. Since  the  
   $\rho_i$'s are
    distinct positive numbers we deduce that for any
     linear map $L:\R^{n-1}\to\R$ we have
 $$\int x_i L(x^{(i)})I_E(U_0(x))f'(|x|_E^2)g(|x|)\d x =0.
 \leqno (2)$$   
 
 For $j\le n$  let $u_j$ be the $j$-th row of $U_0$ and $u_{(j,s)}$ the $s$-th 
 element of $u_j$. For $y\in\R^{n-1}$ we define
  $$\eqalign{ L(y):=&\left(\sum_{j=1}^n u_{(j,i)}^2/r_j^2\right)^{-1}
                       \sum_{j=1}^n {u_{(j,i)}\over r_j^2} <u^{(i)}_j,y>
	                   {\rm\ and}  \cr
	      Q(y):=&\left(\sum_{j=1}^n u_{(j,i)}^2/r_j^2\right)^{-1}		     
                  \left(\sum_{j=1}^n { <u^{(i)}_j,y>^2\over r_j^2}-1\right).
		    \cr}$$ 
  For $x\in\Rn$ we observe that the following equivalences hold.
  $$\eqalign{ &U_0(x)\in E   \cr
      &\iff \sum_{j=1}^n r_j^{-2}[u_{(j,i)} x_i +<u^{(i)}_j,x^{(i)}>]^2\le 1\cr
	 &\iff x_i^2\sum_{j=1}^n u_{(j,i)}^2r_j^{-2} 
	      +2x_i \sum_{j=1}^n u_{(j,i)} r_j^{-2} <u^{(i)}_j,x^{(i)}>
              +\sum_{j=1}^n <u^{(i)}_j,x^{(i)}>^2 r_j^{-2}\le 1 \cr
	&\iff x_i^2+ 2x_iL(x^{(i)})+Q(x^{(i)})\le 0  \cr
	&\iff L^2(x^{(i)})\ge Q(x^{(i)}){\rm\ and\ }
	  |x_i+L(x^{(i)})|\le \sqrt{L^2(x^{(i)})- Q(x^{(i)})}.\cr }$$       
  We claim that $L\equiv0$. Indeed, from the  equivalences above and (2) 
  we deduce  
   that
   $$\eqalign{ 0=& \int_{\{x: U_0(x)\in E\}} 
                  x_i L(x^{(i)}) f'(|x|_F^2) g(|x|) \d x \cr
               =& \int_{L^2(x^{(i)})\ge Q(x^{(i)})} L(x^{(i)})       
 \left[\int_{-L(x^{(i)})-\sqrt{L^2(x^{(i)})-Q(x^{(i)})}}^{-L(x^{(i)})+\sqrt{L^2(x^{(i)})- Q(x^{(i)})}} x_i f'(|x|_F^2) g(|x|) \d x_i\right]\d x^{(i)} .\cr}$$
  Since  for fixed $x^{(i)}$ the function 
   $x_i\mapsto  x_i f'(|x|_F^2) g(|x|)$ is odd and positive if and only if 
   $x_i$ is negative we deduce that  
 $$\int_{-L(x^{(i)})-\sqrt{L^2(x^{(i)})-Q(x^{(i)})}}^{-L(x^{(i)})+\sqrt{L^2(x^{(i)})- Q(x^{(i)})}} x_i f'(|x|_F^2) g(|x|) \d x_i $$
  is negative (respectively, positive) if and only if  
   $L(x^{(i)})$ is positive (respectively, negative). Thus we deduce that 
 $$L(x^{(i)})\int_{-L(x^{(i)})-\sqrt{L^2(x^{(i)})-Q(x^{(i)})}}^{-L(x^{(i)})+\sqrt{L^2(x^{(i)})-Q(x^{(i)})}} x_i f'(|x|_F^2) g(|x|) \d x_i $$
  is negative if and only if $L(x^{(i)})\not= 0$ and vanishes otherwise. Since
   $Q(0)<0$ the inequality ${L^2(x^{(i)})\ge Q(x^{(i)})}$ has solutions for 
    a neighborhood of $0$. This forces $L\equiv0$. 
    Going back to the definition of $L$  we just showed that for
    $\ell\not= i$ the  $\ell$-th coordinate  of 
     $$  \sum_{j=1}^n{u_{(j,i)}\over r_j^2} u_j$$
    vanishes. But, on the other hand this coordinate is equal to the
    element in the $i$-th row and $\ell$-th column of the  product 
     $$U_0^T\circ
 \pmatrix{   r_1^{-2} &   {}  &{}   \cr
              {}      &\ddots  &{}   \cr
	      {}      &   {}  &r_n^{-2}\cr} \circ U_0.$$
 Since $i\not=\ell$ are arbitrary elements of $\{1,\ldots,n\}$ this
  says that above product is a diagonal matrix which finishes the proof of
   the theorem.\bbbox	      
\ms
While we do not know if C$^{\prm}$ holds for an arbitrary $g$ and $f=I_E$, 
where $E$ is an ellipsoid, 
we show below that it does hold for $f$ being a Gaussian density, and $g$  
 log concave.  
\ms
\proclaim Proposition 11. If $g$ is a non-negative, symmetric, log-concave function on $\Rn$ and $A$ is a non-negative definite matrix,  then 
$$\E_\mu \bigl[ \exp(-{\dsp {1\over 2}}<Ax,x>)g(x)\bigr]  \ge \E_\mu \bigl[ \exp(-{\dsp {1\over 2}}<Ax,x>)\bigr]  \E_\mu \bigl[ g(x)\bigr]  .$$

\pf It suffices to assume that $\mu=\mu_n$. Then, 
$$\E_\mu \bigl[ \exp(-{\dsp {1\over 2}}<Ax,x>)g(x)\bigr]  
=(\det(I+A))^{-1/2}\E_\mu \bigl[ g((I+A)^{-1/2}(x))\bigr]  .$$ \par  
We now diagonalize $(I+A)^{-1/2}$  
with the unitary $U$, let $h=g\circ U$ and use the fact that $\mu$ is rotation invariant to allow us to write 
$$\E_\mu \bigl[ g((I+A)^{-1/2}(x))\bigr]  
=\E_\mu \bigl[ g((UU^T(I+A)^{-1/2}UU^T(x))\bigr]  
=\E_\mu \bigl[ h(D(x))\bigr]  .$$
So in order to show that 
$$\E_\mu \bigl[ \exp(-{\dsp {1\over 2}}<Ax,x>)g(x)\bigr]  \ge \E_\mu \bigl[ \exp(-{\dsp {1\over 2}}<Ax,x>)\bigr]  \E_\mu \bigl[ g(x)\bigr],$$
we need only show that 
$$\E_\mu \bigl[ h(D(x))\bigr]  \ge E_\mu \bigl[ h(x)\bigr]  .$$
This is intuitively clear since the eigenvalues of our $D$, say, $\{d_i\}_i$ are between $0$ and $1$, and hence pull $x$ closer to the origin, where $h$ is larger. 
To prove this it suffices, by iteration, to consider the case that all but one of the diagonal entries of $D$ are $1$ and the other, say the $i^{\hbox{th}}$,  is between $0$ and $1$, and, then, to prove that for any non-negative, symmetric log-concave function $\phi$ on $\Rn$ we have 
$$\E_\mu \bigl[ \phi(D(x))\bigr]  \ge E_\mu \bigl[ \phi(x)\bigr] .$$

To do this we integrate out all the variables except the $i^{\hbox{th}}$ on both sides. Then, on each side, we have an integral of a symmetric, log-concave function (by Prekopa-Leindler) of one variable. Such a one dimensional function is non-increasing on the positive axis. Hence, we have a pointwise inequality on the two one-dimensional functions, which yields the inequality we wanted. \bbbox

\bigskip

\centerline{\bf Bibliography}
\bigskip
\item{
[AD]} R.~Ahlswede and D.~E.~Daykin, {\sl An inequality for weight of two
families
 of sets, their unions and intersections,} {\bf Z.
Wahrscheinlichkeitstheorie
 und Verw. Gebiete 43} (1978), 183--185.
\item{[BL]} H.~J.~Brascamp and E.~H.~Lieb, {\sl On the extensions of
 the Brunn-Minkowski and Pr\'ekopa-Leindler theorems, including inequalities
  for log concave functions, and with an application to the diffusion 
   equation}, {\bf J. Functional Analysis 22} (1976), 366--89.
\item{
[DEOPSS]} S.~Das~Gupta, M.~L.~Eaton, I.~Olkin, M.~Perlman, L.~J.~Savage and
 M.~Sobel, {\sl Inequalities on the probability content of convex regions
for
 elliptically contoured distributions,} {\bf proc. Sixth Berkeley Symp.
Math.
 Statist. Prob. 3} (1972), 241--264.

\item{ [D]} O.~J.~Dunn, {\sl Estimation of the means of dependent
variables,} {\bf Ann.
 Math. Statist. 29} (1958), 1095--1111.

\item{ [DS]} C.~W.~Dunnett and M.~Sobel, {\sl Approximations to the
probability integral
 and}\hfil\break {\sl certain percentage points to a multivariate analogue
of
 Student's t-distribution,} \hfil\break {\bf Biometrika 42} (1955),
258--260.

\item{ [Gl]} E.~D.~Gluskin, {\sl Extremal properties of orthogonal
parallelepipeds and
 their applications to the geometry of Banach spaces,} {\bf Math. USSR
Sbornik
 64} (1989), 85--96.

\item{ [KR]} S.~Karlin and Y.~Rinott, {\sl Classes of orderings of measures
and related
 correlation inequalities. I. Multivariate totally positive
distributions,} {\bf
 Journal of Multivariate Analysis 10} (1980), 467--498.

\item{ [Kh]} C.~G.~Khatri, {\sl On certain inequalities for normal
distributions and
 their applications to simultaneous confidence bounds,} {\bf Ann. Math.
Statist.
 38} (1967), 1853--1867.

\item{ [Le]} L.~Leindler, {\sl  On a certain converse of H\"older's
                inequality II,} {\bf Acta. Sci. Math. Szeged 33}, (1972),
 217--223.

\item{ [Pi]} L.~D.~Pitt, {\sl A Gaussian correlation inequality for
               for symmetric convex sets,} {\bf Annals of Probability  5}
               (1977), 470--474.

\item{ [Pr]} A.~Pr\'ekopa, {\sl On logarithmic concave measures and
              functions,} {\bf Acta Sci. Math. Szeged 34} (1973), 335--343.

\item{ [Si1]} Z.~\u Sid\' ak, {\sl Rectangular confidence regions for the
means of
 multivariate normal distributions,} {\bf J. Amer. Statist. Assoc. 62}
(1967),
 626--633.

\item{ [Si2]} Z.~\u Sid\' ak, {\sl On multivariate normal probabilities of
rectangles:
 their dependence on correlations,} {\bf Ann. Math. Statist. 39} (1968),
 1425--1434.
\medskip
{\settabs 3 \columns 
\+G. Schechtman&& Th. Schlumprecht \cr
\+Department of Theoretical Mathematics&&J. Zinn\cr
\+Weizmann Institute of Science&&Department of Mathematics\cr
\+Rehovot, Israel&&Texas A\&M University\cr
\+&&College Station, TX 77843\cr}
\vfill\eject
	     
\bye